\documentclass[a4paper,notitlepage,oneside,openright,12pt]{report}
\usepackage[english]{babel}
\usepackage{latexsym}
\usepackage{amssymb}

\newtheorem{theorem}{Theorem}[section]

\newtheorem{proposition}[theorem]{Proposition}
\newtheorem{remark}[theorem]{Remark}
\newtheorem{definition}[theorem]{Definition}

\renewenvironment{abstract}{\indent\begin{footnotesize}{\bf
Abstract. } }{\end{footnotesize}}
\newenvironment{references}[1]{\begin{footnotesize}
\begin{thebibliography}{#1}}{\end{thebibliography}\end{footnotesize}
}

\begin{document}
\font\titdtr=cmti10 scaled \magstep2
\baselineskip=14pt

\centerline{\sc Families of maximal subbundles of stable
vector bundles on curves}

\vspace{1cm}

{\centerline {\sc Edoardo Ballico and  Barbara Russo}}

\begin{footnote}{AMSC: 14H60, 14D20.
e-mail: russo@degiorgi.science.unitn.it}\end{footnote}

\def\x{$X\,$}
\def\dim{\mbox{dim}}
\def\rk{$ \mbox{rk }$}
\def\kk{\mbox{Ker }}
\def\parn{\par\noindent}
\def\sp{\parn {\bf Step }}
\def\ii{\mbox{Im }}
\def\Z{ Z \!\!\! Z}
\def\Bbb{K \!\!\!\!\! I}
\vspace{2cm}

\begin{abstract}
Let $X$ be a smooth projective curve of genus $g\geq 2$ and let $E$ be a vector bundle on $X.$ Let $M_k(E)$ be the scheme of
all rank $k$ subbundles of $E$ with maximal degree. For every integers $r,$ $k$ and $x$ with $0<k<r$ and either $2k\leq r$
and
$0\leq x\leq (k-1)(r-2k+1)$ or $2k>r$ and $0\leq x\leq (r-k-1)(2k-r+1),$ we construct a rank $r$ stable vector bundles
$E$ such that $M_k(E)$ has an irreducible component of dimension $x.$ Furthermore, if there exists a stable vector bundle
$F$ with small Lange's invariant $s_k(F)$ and with  $M_k(F)$ `spread enough', then $X$ is a multiple covering of a
curve of genus bigger then $2.$
\end{abstract}

\normalsize

\section{Introduction.}
Let $X$ be a smooth projective curve
of genus $g\geq 2$ defined over an algebraically closed field $\Bbb .$
In this paper we study the rank $r $ stable vector bundles, E, on $X$ such that for some integer $k$
with $0<k<r$ $E$ has a `large' family of subbundles with rank $k$ and maximal degree. For positive integers $r,$  $d$ 
let
$M(X; r, d)$ be the moduli space of stable vector bundles on $X$  of rank $r$ and degree $ d. $ It
is well known that  $M(X; r, d)$ is smooth and irreducible. For a positive integer $k$ with $0<k<r,$
let $M_k(E)$ be the set of all rank $k$ subbundles of
$E$ with maximal degree. Being a Quot-scheme, $M_k(E)$ has a natural scheme-structure. For the intent of this paper we
will only need to consider its reduced structure. Indeed we are interested in finding a stable vector bundle 
$E$ such that
$M_k(E)$ has an irreducible component with prescribed dimension. Since every element in $M_k(E)$ has
maximal degree, the scheme $M_k(E)$ is complete. Hence by \cite{MS}, pp. 254-255, we have
$\dim(M_k(E))\leq k(r-k)$ for every rank $r$ vector bundle $E.$ Fixed $x$ with
$x\leq k(r-k),$ it is very easy to find a decomposable rank $r$ vector bundle $E$ such that $M_k(E)$ has an irreducible
component of dimension $x.$ But  we are interested in stable vector bundles which are indecomposable. Hence using extensions
of a line bundle by a decomposable rank $r-1$ bundle  we will prove in section $2$ the following result:
\begin{theorem}
Fix integers $g,$ $r,$ $k$ with $2\leq g\leq r+1,$ $0<k<r$; if $2k\leq r,$ then assume $x\leq (k-1)(r-2k+1);$ if $2k>r,$
then assume $x\leq (r-k-1)(2k-r+1).$ Let $X$ be a smooth projective curve of genus $g.$ Then there exists a stable
vector bundle $E$ on $X$ such that $M_k(E) $ has an irreducible component of dimension $x.$ \label{teorema}   
\end{theorem}
\vspace{0,5cm}\par

The proof of Theorem \ref{teorema} is quite simple but even if we tried we were not able to produce larger families of maximal
degree subbundles. The bound on the dimension $x:=\dim(M_k(E))$ seems to be quite good, (see Proposition \ref{teorema3}). The
dimension of
$M_k(E)$ is known when
$E$ is a general element of
$M(X;r,d)$  (see Remark \ref{remark} and
Proposition \ref{proposizione}). Classically   the picture was clear
for a rank $2$ stable vector bundle $E:$ either $\dim( M_1(E))=0$ or    $\dim( M_1(E))=1 $ (see the introduction of \cite{LN}
and references therein). In fact 
the situation is described
by one invariant, called degree of stability, $s(E).$  It is known that $0<s(E)\leq g$ and  
$s(E)\simeq \deg(E) \ (2)$
(\cite{Na}).  Furthermore, for $E$ general in its moduli space we have
$s(E)=g$ if
$g-d$ is even and $s(E)=g-1$ if $g-d$ is odd. Maruyama proved two main facts:  if $s(E)=g,$  then
$\dim( M_1(E))=1$ and  if
$s(E)<g
$ then dim$(M_1(E))=0$. H. Lange and M.S. Narasimhan
produced  examples of stable rank 2 vector bundles with
$\dim( M_1(E))=0$ and
$s(E)<g$ 
 (see \cite{LN}, Prop. 3.3. and sections $5,$ $6$ and $7).$  Indeed taking $f:X\to Y$ a multiple
covering of curve
$Y$  of genus
$g'\geq 2$ they  were able to produce  examples of curves $X$ of genus $g$  big enough to obtain a   stable rank
$2$ vector bundle, $E,$ on $X$ with
$s(E)<g$ and
$\dim( M_1(E))=1,$ by pulling back a stable vector bundle, $F,$ on $Y$ with $s(F)=g'$ 
(see \cite{LN}, Prop. 7.3).  In \cite{Bu} D. Butler proved some kind of reverse question: 
  if
$E$ is a stable vector bundle of rank $2$ with 
$\dim(M_1(E))
=1$ and
$s(E)(2s(E)-1)<g$ then there is a covering $f:X\to Y$ and a stable vector bundle on $Y,$ $F,$  
with $R\in \mbox{Pic}(X)$ with
$A\otimes R\simeq f^* (B)$  and $\dim(M_1(F))=1. $ In higher rank the situation is more
complicated (see Remark \ref{remark}). In particular  the stability condition for a rank $r$
vector bundle,
$E,$   is controlled by
$r-1$ invariants called degrees of stability (or Lange's invariants): $$s_k(E)=k\deg
(E)-r\min_{\scriptstyle \begin{array}{l}H\hookrightarrow E\\  \rk H=k\end{array}}\deg(H).$$ In section
$3$ we give a partial generalization to higher rank of a theorem of D. C. Butler (see Theorem
\ref{teorema2}) which gives how restrictive is to have `many and very spread' maximal degree
subbundles. This is the key motivation of our paper: Theorem \ref{teorema2} and Proposition \ref{teorema3} show the 
existence of a rank
$r$ stable vector bundle,
$E,$ with a low value of $s_k(E)$ and large dimension of $M_k(E).$
\vspace{0.5cm}

This research was partially supported by MURST. Both authors are members of the VBAC
Research group of Europroj.

\section{Proof of Theorem \ref{teorema}}

Before proving Theorem \ref{teorema} we need the following remark
\begin{remark} Assume char$\Bbb  =0.$ Fix some integers $g,$ $r,$ $k,$ $a,$ $b$ 
with $g\geq 3,$ $r\geq 2,$ $0<k<r$ and
$kb-a(r-k)>0.$ Let $X$ be  a smooth projective curves of genus $g.$ Let $A$ be  
a general member of $M(X;k,a),$    $B$ 
a general member of $M(X;r-k,b)$ and $E$ a general extension of $B$ by $A.$ If 
$kb-a(r-k)<k(r-k)(g-1),$ by
\cite{RT}, Thm.0.1,
$E$ is stable (see also 
\cite{BL} for several special cases).  Furthermore,  by a result of 
A. Hirschovitz (\cite{Hi}) a general member
of $M(X;r,a+b)$ is an extension of a general  $B\in M(X;r-k,b)$ by a general  $A \in M(X;k,a)$ if and only if $kb-a(r-k)\geq
k(r-k)(g-1).$
   As remarked in the
introductions of \cite{RT} and \cite{BL} (\cite{BL} eq. (D)), the stability of 
such an $E$ implies 
dim$(M_k(E))=\max\{ s-k(r-k)(g-1), 0\}. $ In fact $M_k(E) $ turns out to be the 
fiber of a morphism,
$\phi,$ between the parameter space of stable extensions of stable vector bundles and the moduli space $M(X;r,d);$   this
allows to estimate the dimension of $M_k(E).$  In particular if
$s=k(r-k)g$ then 
$\dim( M_k(E))=k(r-k)$ which by
\cite{MS}, pp 254-255, it is the maximum admissible dimension of  $M_k(E).$ \label{remark}
\end{remark}
\vspace{0,5cm}\par
If char$\Bbb  =0$ there exists a first weak version of theorem \ref{teorema}:
\begin{proposition}
Assume char$\Bbb  =0$. Fix integers $r,$ $k,$ $x$ with $0<k<r,$  $0\leq x\leq k(r-k)$ and $x$ 
divisible by the highest common
divisor, $u,$ of $k$ and $r.$ Let $X$ be a smooth curve of genus $g\geq 3.$ Then there exists an
integer $d$ such that for a general $E\in M(X;d,r)$ the algebraic set $M_k(E)$ has an irreducible
component of dimension $x$ and every irreducible component of $M_k(E)$ ha dimension at most $x.$
\label{proposizione}
\end{proposition}
\vspace{0,5cm}\par
{\it Proof.}  Since $u$ divides $x,$ there exists an integer
$d$ with
$0\leq d<r.$ Moreover there exists an unique integer $a$ satisfying $\frac{d-a}{r-k}-g\leq \frac{a}{k}\leq
\frac{d-a}{r-k}-g+1.$ Hence as pointed out in \ref{remark} we have $\dim (M_k(E))=x=\max\{ s-k(r-k)(g-1), 0\}$ with
$s=(d-a)k-a(r-k).$ 
\vspace{0,5cm}\par
{\it Proof of \ref{teorema}}  Since the cases $k=1$ and $k=r-1$ are covered by Proposition \ref{proposizione}, when
char$\Bbb  =0$ and $g\geq 3,$ we may assume $k\geq 2$ and $r-k\geq 2.$ Furthermore, $M_k(E)\simeq M_{r-k}(E^*)$ for every
rank
$r$ vector bundle $E.$ Therefore taking, if necessary, the dual bundle, we may assume $2k\leq r.$ 
If char$\Bbb  >0$ or $g=2$ and
$k=1$ or $k=r-1$ proceed as in the last part of case 2) below.   
Hence from now on we may assume $4\leq 2k\leq r.$ Since
$x\leq (k-1)((r-k)-(k-1))$ we can find two integers $y$ and $t$ with $0<2t\leq y\leq r-k,$ 
$t\leq k-1$ and
$t(y-1-t)\leq x \leq t(y-t).$ Set $e:=x-t(y-1-t)$. Then $0\leq e<t$ and if $y=r-k$ then $e=0.$  
Therefore $y+e+1\leq r-1.$ Take a general $(r-e-y-1)$-ple
$(M,R_1,...,R_{r-e-y-1})\in \mbox{Pic}^0 (X)\times  ...
\times 
\mbox{Pic}^0 (X)$ and $L\in  \mbox{Pic}^1 (X)$ with $h^0(X,L)=0.$
Set $F:= {\cal O}_X^{\oplus y}\oplus M^{\oplus (e+1)}\oplus (\oplus _{1\leq i\leq r-e-y-1} R_i)$ (notice that $y+e+1\leq r-1)$.
By construction $F$ is a semi-stable vector bundle with $\rk F=r-1$ and $\deg F=0.$ Let $E$ be a general extension of $L$
by $F.$ 

{\bf Claim.}  $E$ has no proper subsheaf with positive degree and every degree $0$ subsheaf
of $E$ is a subsheaf of
$F.$

Here we assume the Claim. Hence $E$ is stable. Choose some integers $u,$ $v$ with $0\leq u \leq y,$ $0\leq v\leq e+1$
and 
$0\leq k-u-v\leq r-e-y-2.$  Let $I$ any subset of $\{1,..., r-e-y-2\}$ with card$(I)=k-u-v.$ Call
$T(u,v,I)$ the following family of rank $k$ subbundles of $F$ with degree $0:$ $A\in T(u,v,I)$ if
and only if $A\simeq A_1\oplus A_2\oplus A_3$ where $A_1$  subsheaf of ${\cal
O}^{\oplus y}$ isomorphic to ${\cal
O}^{\oplus u}$  $A_2$ is a subsheaf of $M^{\oplus (e+1)}$ isomorphic to $M^{v}$ and
$A_3\simeq
\oplus _{i\in I} R_i.$  Since $F$ is polystable and no two among the degree $0$ line bundles
${\cal O}_X,$ $M$ and $R_i,$ $1\leq i\leq r-y-e-2,$ are isomorphic, then $T(u,v,I)$ is an irreducible component of
$M_k(E)$ with $\dim  ( T(u,v,I))=u(y-u)+(e+1-v)v.$ Varying $u,$ $v$ and $I$ we obtain in this way all
the irreducible components of $M_k(F).$ 
By the second part of the Claim, these are the irreducible components of $M_k(E).$ When $u=t $ and 
$ v=1$ by the definition of $e$ we get $\dim(T(t,1,I))=x.$  Hence to prove \ref{teorema} it is
sufficient to prove the Claim.
\vspace{0,2cm}\par
{\it Proof of the Claim.} We move the line bundles $M$ and $R_i$ $1\leq i\leq r-e-y-2,$ in
Pic$^0(X).$ By the semicontinuity of the Lange's invariants $s_k$ (\cite{La}, Lemma 1.3), it is
sufficient to prove the Claim  for the following  general extension,
\begin{eqnarray}0\to {\cal O}_X^{\oplus (r-1)} \to G\to L\to 0 \label{equazione}\end{eqnarray}.

Since $h^0(X,L) =0,$ we
have
$h^0(X,G)=r-1.$ In particular the subsheaf ${\cal O}_X^{\oplus (r-1)} $ is the subsheaf spanned by
$H^0(X,G).$ Hence it is uniquely determinated by $G$ and sent into itself by any endomorphism of
$G.$ Therefore $G$ fits in a unique way into \ref{equazione}, up to an element of Aut$(G).$ Since
$\chi (L^*) =-g$ and by our assumptions on $g$ and $r,$ $G$  contains no factor isomorphic to
${\cal O}_X.$ In order to obtain a contradiction we assume
the existence of a proper subsheaf
$B$ of
$G$ with $\deg (B)\geq 0,$ and if $\deg B =0$ we suppose that $B$ is not a direct factor of ${\cal
O}_X ^{\oplus (r-1)}.$ Taking $h:=\rk B$ minimal among all the ranks of such subbundles, we may assume
$B$ stable. Taking $\deg (B)$ maximum among all the degrees of all such rank $h$ subbundles we may
assume $B$ saturated in $G.$ Since $B$ is not contained in  ${\cal
O}_X ^{\oplus (r-1)},$ the map $\pi : B \to L $ induced by the surjection $j:G\to L$ in
\ref{equazione} is not zero. Set $B':  \kk (\pi),$  $L'=\ii (\pi)$ and $w:=h^0(X,B').$
Since $B'$ is a subsheaf of  ${\cal
O}_X ^{\oplus (r-1)},$ we have $B'\simeq B''\oplus  {\cal O}_X ^{\oplus (w)} $ with $h^0(X,B'')=0.$
Since $B'^{*}$ is spanned, $\det (A'^*)$ is spanned. Thus if $\deg (B'^*) =\deg (\det
(B'^* )\not =0,$ $X$ has a degree $\deg(B'^*)$ pencil. By our assumption on the degree of $B$ we
have  $\deg (B'^*)\leq \deg (L') \leq \deg (L) =1. $ Since $g>0$ there is no degree $\deg (B'^*)$
pencil on $X.$ Hence a contradiction.  Thus
$\deg(B'^*)=0,$ that is $w=h-1$ and $B'\simeq {\cal O}_X ^{\oplus (h-1)}.$ At this point we
distinguish two cases:

{\it Case 1)} Here we assume $L\not \simeq L',$ that is the existence of a positive divisor $D$ with
$L'=L(-D).$ Since $\deg (L')\leq \deg (L)-1=0,$ $\mu (B)\geq 0$ and $B$ is stable, we obtain
a contradiction, unless $h=1,$ $B\simeq L'$ and $w=0.$ In this case we have $L'\simeq L(-P)$ for some
$P\in W$ and $F$ a positive  elementary transformation of  ${\cal O}_X ^{\oplus (r-1)}\oplus L(-P)$
supported in $P.$ Hence the set of all such bundles $G$ depends at most on $r$ parameters. Since
$\dim(\mbox{Ext}^1(L,{\cal O}_X ^{\oplus (r-1)})))=(r-1)g$ by the Riemann-Roch Theorem and any such
$G$ fits, up to a multiplicative constant, in a unique exact sequence \ref{equazione}, we get a
contradiction concluding the proof in Case 1).

{\it Case 2)} Here we assume $L\simeq L'.$ Then since $B'\simeq {\cal O}_X ^{\oplus (w)}$ as direct
factor of ${\cal O}_X ^{\oplus (r-1)}$  we get
$G/B\simeq {\cal O}_X ^{\oplus (r-1-w)}={\cal O}_X ^{\oplus (r-h)}.$ Hence $G/B$ is isomorphic to a
direct factor of $G.$ But $G$ cannot have any trivial  factor which is a contradiction and the
theorem is proved. 

\begin{remark} The proof of \ref{teorema} shows the existence of a vector bundle $E\in M(X;r,1)$ such
that $M_k(E)$ has an irreducible component $t$ of dimension $x$ and such that every $B\in T$ is a
direct sum of line bundles of degree $0.$
\end{remark}
\vspace{0,2cm}\par
\begin{remark}
Let $T\subset M_k(E)$ be an irreducible subvariety such that there is a subbundle $F$ of $E$
containing every $B\in T.$ By \cite{MS}, pp. 254-255, it follows $\dim(T)\leq k(r-k).$ In the proof
of Theorem \ref{teorema} we have constructed a vector bundle $E$ which has a subbundle $F$ with
exactly this property.  
\end{remark}
\vspace{0,5cm}\par
We repeat here the description of the irreducible components of $M_k(E)$ for the stable bundle, $E,$
obtained in the proof of Theorem \ref{teorema}. First choose integers $u,$ $v$ with $0\leq u\leq y,$
$0\leq v\leq e+1,$ $0\leq k-u-v\leq r-e-y-2.$  Then choose any subset, $I,$ of $\{ 1, ...,
r-e-y-2\}$ with card$(I)=k-u-v.$ For any such data $(u,v,I)$ there is an irreducible component,
$T(u,v,I)$ of $M_k(E)$ and every irreducible component of $M_k(E)$ arises in this way. Furthermore, we
have $\dim (T(u,v,I))=u(y-u)+(e+1-v)v.$
\section{Maximally spread families and multiple covering curves}

In this section we will give a partial generalization of a result of D. C. Butler, \cite{Bu}.
As in \cite{Bu} we will use a result of Accola (\cite{Ac}) which is valid in characteristic zero.
Therefore we assume char$\Bbb  =0.$ Let $X$ be a smooth projective curve of genus $g\geq 2.$ Fix
two integers $k,$ $r$ with $0<k<r$ and set $m:=GCD(k,r-k),$ $v:=\frac{r-k}{m}$ and $w:=\frac{k}{m}.$
Let $E$ be a rank $r$ vector bundle on $X$ and ${\cal H}:=\{H_t\}_{t\in T}$ be a flat family of
saturated rank $k$ subbundles of $E$ parameterized by an irreducible complete variety $T.$ For every
$t\in T$ set $G_t:= E/H_t.$ For all pairs $(x,y)\in T^2$ the composition of the inclusion
$i_x:H_x\to E$ with the surjection $j_x: E\to G_y$ gives a map $\phi(x,y):H_x\to G_y$ such that
$\phi(x,y)=0 $ if and only if $H_x$ and $H_y$ are isomorphic subsheaf of $E.$ More generally, for
all $(x(1),...,x(v),y(1),...,y(w))\in T^{v+w}$ we have a map $\Phi((x(1),...,x(v),y(1),...,y(w))):
H_{x(1)}\oplus ...\oplus H_{ x(v)}\to G_{y(1)}\oplus ...\oplus G_{y(w)}
.$  Notice that $H_{x(1)}\oplus ...\oplus H_{ x(v)}$ and $ G_{y(1)}\oplus ...\oplus G_{y(w)}$ have
the same rank $\frac{k(r-k)}{m}.$ 
\begin{definition} The family $\cal H$ is called {\it maximal spread} if for general
$(x(1),...,x(v),y(1),...,y(w))\in T^{v+w}$ the map $\Phi((x(1),...,x(v),y(1),...,y(w)))$ is
invertible at a general point of $X.$
\end{definition}
\begin{remark} If $r=2k$ maximally spread means that for general $(x(1),y(1))\in T^{2}$ the map
$H_{x(1)}\to G_{y(1)} $ is an injective map of sheaves, which is a condition that may be satisfied.
\end{remark}\vspace{0,5cm}\par
By definition a maximal spread family $\cal H$ induces an inclusion of sheaves of
$H_{x(1)}\oplus ...\oplus H_{ x(v)}$ in $ G_{y(1)}\oplus ...\oplus G_{y(w)}.$ 
If $\cal H$ is maximal spread then the map 

\begin{displaymath} \det(\Phi((x(1),...,x(v),y(1),...,y(w)))):
\det(H_{x(1)}\oplus ...\oplus H_{ x(v)})\to \det(G_{y(1)}\oplus ...\oplus G_{y(w)})
\end{displaymath} 

is an inclusion. Therefore there is an effective divisor, $Z((x(1),...,x(v),y(1),...,y(w))),$
associated to a line bundle isomorphic to $\det(H_{x(1)}\oplus ...\oplus H_{ x(v)})^*\otimes
\det(G_{y(1)}\oplus ...\oplus G_{y(w)}).$ Hence
$$\begin{array}{l}\deg(Z((x(1),...,x(v),y(1),...,y(w))))=w(\deg(G_t))-v(\deg(H_t))=
\\ \\ =w(\deg(E)-\deg(H_t))-
v(\deg(H_t)=\frac{(k(\deg(E)-r(\deg(H_t))}{m}.\end{array}$$
Hence if $H_t$ is maximal (that is has maximum degree among rank $k$ subbundles of $E)$  then 
$\deg(Z((x(1),...,x(v),y(1),...,y(w))))=\frac{s_k(E)}{m}.$ 
The divisor
$Z((x(1),...,x(v),y(1),...,y(w)))$ depends symmetrically on the variables $x(i)\in T,$ $1\leq i\leq
v,$ and $y(j)\in T,$ $1\leq j\leq
w.$
Notice that we have defined the divisors $Z((x(1),...,x(v),y(1),...,y(w)))$ in a general open set of
$T^{v+w}.$ Since $T$ is complete the set of effective divisors 
$Z((x(1),...,x(v),y(1),...,y(w)))$ has
limits for all $(x(1),...,x(v),y(1),...,y(w))\in T^{v+w}.$ These limits are not unique, but this does
not effect  our computation. In particular for every  $x\in T,$ we may find $Z(x,...,x,x,...,x)$
an effective divisor such that 
${\cal O}(Z(x,...,x,x,...,x))\simeq \det(H_x)^{\otimes v}\otimes \det(G_x)^{\otimes w}.$

\begin{remark}
Notice that for every $(x(1),...,x(v),y(1),...,y(w))\in T^{v+w}$ the divisor $$(v+w)Z((x(1),...,x(v),y(1),...,y(w)))$$ and the
divisor $$\sum_{1\leq i\leq v}Z((x(i),...,x(i),x(i),...,x(i)))+\sum_{0\leq j \leq w}Z((y(j),...,y(j),y(j),...,y(j)))$$ are
associated to the same line bundle $$\det (H_{x(1)}\oplus ...\oplus H_{ x(v)})^*\otimes
\det(G_{y(1)}\oplus ...\oplus G_{y(w)})^{(v+w)}$$ and therefore they are linearly equivalent.
Call $L((x(1),...,x(v),y(1),...,y(w)))$ the subsheaf of $\det(H_{x(1)}\oplus ...\oplus H_{ x(v)})^*\otimes
\det(G_{y(1)}\oplus ...\oplus G_{y(w)})^{(v+w)}$ spanned by $H^0(X, \det(H_{x(1)}\oplus ...\oplus H_{
x(v)})^*\otimes
\det(G_{y(1)}\oplus ...\oplus G_{y(w)})^{(v+w)}).$ We believe that the two  families of line bundles
$\{\det(H_{x(1)}\oplus ...\oplus H_{ x(v)})^*\otimes
\det(G_{y(1)}\oplus ...\oplus G_{y(w)})\}$ and $\{L((x(1),...,x(v),y(1),...,y(w)))\ | \
(x(1),...,x(v),y(1),...,y(w))\in T^{v+w}\} $ give more information on the geometry of $E$ then
$s_k(E)$ (even in the case in which $M_k(E)$ is finite).\label{remark1}
\end{remark}
\begin{theorem}
Assume char$\Bbb   =0.$ Let $X$ be a smooth projective curve of genus $g\geq 2$ and $E\in
M(X;r,d),$ $r\geq 2,$ such that $M_k(E)$ has a maximal spread family,  $T,$ and such that
$s_k(E)(  s_k(E)-m)<m^2g,$ where $m:=GCD(k,r).$ Then there exist a smooth curve $C$ and a morphism
$\pi:X\to C$ with $\deg(\pi)>1.$
\label{teorema2}\end{theorem}
\begin{remark} As one can easily see we are going to prove more then 
 what is stated in the Theorem \ref{teorema2}. In fact we are going to prove that there exists a family of line bundles 
$R(x(1),...,x(v),y(1),...,y(w))\in \mbox{ Pic }(C)$ such that $\pi^*(R(x(1),...,x(v),y(1),...,y(w)))\simeq 
\det(H_{x(1)}\oplus ...\oplus H_{ x(v)})^*\otimes
\det(G_{y(1)}\oplus ...\oplus G_{y(w)}).$ If the rank of $E$ is $2$ the existence of this family (with $w=v=1)$ allows to
construct a rank
$2$ stable vector bundle $F$ on $C$ whose pull-back is $E$ and whose family of maximal degree linebundles is the pull-back of
the one of $E,$  up to a twist by a line bundle, $A,$ on $C,$ (see \cite{Bu}). 
\end{remark}
\vspace{0,5cm}\par
{\it Proof.} Set $v:=\frac{r-k}{m}$ and $w:=\frac{k}{m}$ and take general
$(x(1),...,x(v),y(1),...,y(w))\in T^{v+w}.$  By Remark \ref{remark1} we have 
$$h^0(\det(H_{x(1)}\oplus ...\oplus H_{x(v)})^*\otimes
\det(G_{y(1)}\oplus ...\oplus G_{y(w)})^{(v+w)})\geq 2 .$$ As in Remark \ref{remark1} consider the
line bundles $L((x(1),...,x(v),y(1),...,y(w)));$ they form an infinite family of spanned non-trivial
line bundles with degreee at most $\frac{s_k(E)}{m}.$ Since $\frac{s_k(E)}{m}(\frac{s_k(E)}{m}-1)<g,$ we
can apply a result of Accola (see
\cite{Ac}, Th. 4.3, or \cite{Bu}, Lemma 1.2)  
finding a non-trivial covering $\pi:X\to C$ and
$R(x(1),...,x(v),y(1),...,y(w))\in \mbox{ Pic }(C)$ with 
$\pi^*(R(x(1),...,x(v),y(1),...,y(w)))\simeq 
\det(H_{x(1)}\oplus ...\oplus H_{ x(v)})^*\otimes
\det(G_{y(1)}\oplus ...\oplus G_{y(w)}).$
\vspace{0,5cm}\par 
To explain the notion of maximally spread family, we prove the following easy result
\begin{proposition}
(any char$\Bbb  )$ Let $X$ be a smooth projective curve of genus $g\geq 2.$ Fis integers $r,$ $k$
with $0<k<r$ and a rank $r$ vector bundle $E$ on $X.$ Let $T\subset M_k(E)$ be an irreducible
projective family with dim$(T)>k(r-1-k).$ Then $T$ is maximally spread. Furthermore, for every
$P\in X$ the union of the subspaces $H_t {_{ |_{\{P\}}}}\subset E_{ |_{\{P\}}}$ is not contained in a
lower dimensional vector subspace  of $ E_{ |_{\{P\}}}.$\label{teorema3}
\end{proposition}\vspace{0,5cm}\par
{\it Proof.} Fix $P\in X.$ By the proof of Proposition of pg 254 in \cite{MS}, the map $$\pi:T\to \mbox{Grass
}(r-k,E_{ |_{\{P\}}})$$ sending $H_t,$ $t\in T,$ into the $(r-k)-$dimensionl vector space $E_{
|_{\{P\}}}/H_t{_{ |_{\{P\}}}}$ is finite. Since $\dim(T)>k(r-k)=\mbox{Grass }(r-k,E_
{|\ \{P\}})$ the union of all subspaces $H_t{_{ |_{\{P\}}}}$ for $t\in T$ cannot be contained in a
hyperplane of $E_{ |_{\{P\}}}.$

\par
\vspace{1cm}
E. Ballico, Universit\`a di Trento
\par
Dip. di Matematica - 38050 Povo (TN) - Italy
\par
e-mail: ballico@alpha.science.unitn.it
\par
\vspace{1cm}
E. Russo, Universit\`a di Trento
\par
Dip. di Matematica - 38050 Povo (TN) - Italy
\par
e-mail: russo@degiorgi.science.unitn.it

\end{document}